\documentclass[reqno,a4paper]{amsart}
\usepackage{amsfonts,amssymb,amstext,amsmath,amsthm}

\begin{document}

\newtheorem{thm}{Theorem}
\newtheorem{lem}[thm]{Lemma}
\newtheorem{rem}[thm]{Remark}
\newtheorem{prop}[thm]{Proposition}
\newtheorem{cor}[thm]{Corollary}
\newtheorem{ex}[thm]{Example}

\newcommand{\ri}{r.~i.}


\title[]{Khintchine inequality and Banach-Saks type properties in rearrangement-invariant spaces}
\date{}

\author[]{F.~A.~Sukochev and D. Zanin}

\thanks{Research supported by 
the Australian Research Council.}


\keywords{$p$-Banach Saks type properties, rearrangement-invariant
spaces, Khintchine inequality, Kruglov property}

\subjclass{46E30, 46B20}

\date{}

\begin{abstract}
{\it We study the class of all rearrangement-invariant (=r.i.)
 function spaces $E$ on $[0,1]$ such that there exists
$0<q<1$ for which $ \Vert \sum_{_{k=1}}^n\xi_k\Vert _{E}\leq
Cn^{q}$, where $\{\xi_k\}_{k\ge 1}\subset E$ is an arbitrary
sequence of independent identically distributed symmetric random
variables on $[0,1]$ and $C>0$ does not depend on $n$. We
completely characterize all Lorentz spaces having this property
and complement classical results of Rodin and Semenov for Orlicz
spaces $exp(L_p)$, $p\ge 1$. We further apply our results to the
study of Banach-Saks index
sets in r.i. spaces. 
}
\end{abstract}

\maketitle

\section{Introduction}
A classical result of Rodin and Semenov (see \cite{RS} or \cite[
Theorem 2.b.4]{LT-II}) says that the sequence of Rademacher
functions  $\{r_k\}_{k\ge 1}$ on $[0,1]$ in a r.i. space $E$ is
equivalent to the unit vector basis of $l_2$ if and only if $E$
contains (the separable part of) the Orlicz space $L_{N_2}(0,1)$
(customarily denoted as $exp(L_2)$) where $N_2(t)=e^{t^2}-1$.
Here, $\{r_k\}_{k\ge 1}$ may be thought of as a sequence of
independent identically distributed centered Bernoulli variables
on $[0,1]$. A quick analysis of the proof (see e.g.
\cite[p.134]{LT-II}) shows that the embedding $exp(L_2)\subseteq
E$ is established there under a weaker assumption that
$\{r_k\}_{k\ge 1}$ is $2$-Banach-Saks sequence in $E$, that is $
\Vert \sum_{_{k=1}}^nr_k\Vert _{_{E}}\leq Cn^{1/2}$, where $C>0$
does not depend on $n\ge 1$. The main object of study in the
present article is the class of all r.i. spaces $E$ such that
there exists $0<q<1$ for which
\begin{equation}\label{mainzero} \Vert \sum_{_{k=1}}^n\xi_k\Vert _{E}\leq Cn^{q},
\end{equation}
where $\{\xi_k\}_{k\ge 1}\subset E$ is an arbitrary sequence of
independent identically distributed symmetric random variables on
$[0,1]$ and $C>0$ does not depend on $n$. We completely
characterize all Lorentz spaces from this class in Corollary
\ref{lorentz alternative} below. In Theorem \ref{Marc} we obtain
sharp estimates of type \eqref{mainzero} for the Orlicz spaces
$exp(L_p)=L_{N_p}(0,1)$, $1\leq p<\infty$ where $N_p(t)=e^{t^p}-1$
complementing results of \cite{RS} (see also exposition in
\cite{D}). Our results have also a number of interesting
implications to the study of Banach-Saks type properties in r.i.
spaces.

Recall that a bounded sequence $\{x_n\}\subset E$ is called a
p-BS-sequence if for all subsequences $\{y_k\}\subset\{x_n\}$ we
have
\[
\sup\limits_{m\in
 N}m^{-\frac{1}{p}}\Big\|\sum\limits_{k=1}^{m}y_k\Big\|_E<\infty.
\]
We say that $E$ has the p-BS-property  and we write $E\in BS(p)$
if each weakly null sequence contains a p-BS-sequence.  The set
\[
 \Gamma(E)=\{p:\:p\geq 1,\:E\in BS(p)\}
\]
is said to be the index set of $E$, and is of the form
$[1,\gamma]$, or  $[1,\gamma)$ for some $1\leq \gamma$.
%
%

If, in the preceding definition, we replace all weakly null
sequences by weakly null sequences of independent random variables
(respectively, by weakly null sequences of pairwise disjoint
elements; by weakly null sequences of independent identically
distributed random variables), we obtain the set $\Gamma_{\rm
i}(E)$ (respectively, $\Gamma_{\rm d}(E)$, $\Gamma_{\rm iid}(E)$).
The general problem of describing and comparing the sets
$\Gamma(E)$, $\Gamma_{\rm i}(E)$, $\Gamma_{\rm iid}(E)$) and
$\Gamma_{\rm d}(E)$ in various classes of r.i. spaces was
addressed in \cite{SeSu-CR, DoSeSu2004, SeSu, AsSeSu2005,
new-16-Sem-Suk, AsSeSu2007}. In particular, it is known
\cite{AsSeSu2005} that $1\in\Gamma(E)\subseteq \Gamma_{\rm
i}(E)\subseteq \Gamma_{\rm iid}(E)\subseteq [1,2]$ and
$\Gamma_{\rm i}(E)\subseteq \Gamma_d(E)$ for any r.~i.\ space $E$.
Moreover, the sets $\Gamma(E)$ and $\Gamma_{\rm i}(E)$ coincide in
many cases but not always. For example, $\Gamma(L_p)=\Gamma_{\rm
i}(L_p)=\Gamma_{\rm iid}(L_p)$, $1<p<\infty$ (see e.g.
\cite[Corollary 4.4 and Theorem 4.5]{new-16-Sem-Suk} and also
Theorem \ref{firstmain} below), whereas for the Lorentz space
$L_{2,1}$ generated by the function $\psi(t)=t^{1/2},$ we have
$\Gamma(L_{2,1})=[1,2)$ and $\Gamma_{\rm i}(L_{2,1})=[1,2]$
(\cite[Theorem~5.9]{new-16-Sem-Suk} and
\cite[Proposition~4.12]{AsSeSu2005}). It turns out that these two
situations are typical \cite[Theorem 9]{SeSu}: under the
assumption that $\Gamma(E)\ne \{1\}$, we have either $\Gamma_{\rm
i}(E)\setminus\Gamma(E)=\emptyset$ or else $\Gamma_{\rm
i}(E)\setminus\Gamma(E)=\{2\}$.

The present paper may also be considered as a contribution to the
study of the class of all r.i. spaces $E$ such that $\Gamma_{\rm
iid}(E)=\Gamma_{\rm i}(E)$. We prove a general theorem (see
Theorem \ref{firstmain} below) that $\Gamma_{\rm
iid}(E)=\Gamma_{\rm i}(E)$ if and only if $\Gamma_{\rm
iid}(E)\subseteq \Gamma_{\rm d}(E)$. It is easy to see that every
Lorentz space $\Lambda(\psi)$ satisfies the latter condition and,
using the main result described above, we give a complete
characterization of all Lorentz spaces $E=\Lambda(\psi)$ such that
$\Gamma_{\rm iid}(E)\neq \{1\}$ (see Theorem \ref{mainsecond} and
Corollary \ref{mainsecond_add}).

It also pertinent to note here, that if one views the Rademacher
system as a special example of sequences of independent mean zero
random variables, then a significant generalization of Khintchine
inequality is due to W.B. Johnson and G. Schechtman
\cite{JoSch1989}. They introduced the r.i. space $Z_E^2$ on
$[0,\infty)$ linked with a given r.i. space $E$ on $[0,1]$ and
showed that any sequence $\{f_k\}_{k=1}^\infty$ of independent
mean zero random variables in $E$ is equivalent to the sequence of
its disjoint translates $\{\bar
f_k(\cdot):=f_k(\cdot-k+1)\}_{k=1}^\infty$ in $Z_E^2$, provided
that $E$ contains an $L_p$-space for some $p<\infty$. This study
was taken further in \cite{Br1994,
AsSeSu2005,AsSu2006-1,AsSu2006-2}, where the connection between
this (generalized) Khintchine inequality and the so-called Kruglov
property was established (we explain the latter property in the
next section). We show the connection between the class of all
r.i. spaces with Kruglov property and the estimates
\eqref{mainzero} in Theorem \ref{Kruglov}. Recently, examples of
r.i. spaces $E$ such that $\Gamma(E)= \{1\}$ but $\Gamma_{\rm
i}(E)\ne \{1\}$ have been produced in \cite{AsSeSu2007} under the
assumption that $E$ has the  Kruglov property. Our approach in
this paper complements that of \cite{AsSeSu2007}; in particular,
we present examples of Lorentz and Marcinkiewicz spaces $E$ such
that $\Gamma_{\rm i}(E)=\Gamma_{\rm iid}(E)\neq \{1\}$ and which
do not possess the Kruglov property.


Finally, we show that the equality $\Gamma_{\rm
iid}(E)=\Gamma_{\rm i}(E)$ fails when $E$ is a classical space
$L_{pq}$, $1<q<p<2$.

\section{Definitions and preliminaries}
\subsection{Rearrangement-invariant spaces}
A Banach space $(E,\Vert \cdot\Vert _{_{E}})$ of real-valued
Lebesgue measurable functions (with identification $m$-a.e.) on
the interval $[0,1]$ will be called {\it rearrangement-invariant}
(briefly, \ri ) if
\begin {enumerate}
\item[(i).]  $E$ is an ideal lattice, that is, if $y\in E$, and if
$x$ is any measurable function on $[0,1]$ with $0\leq \vert x\vert
\leq \vert y\vert $ then $x\in E$ and $\Vert x\Vert _{_{E}}
 \leq \Vert y\Vert _{_{E}};$
\item[(ii).]  $E$ is rearrangement invariant in the sense that if
 $y\in E$, and if $x$
is any measurable function on $[0,1]$ with $x^*=y^*$, then $x\in
E$ and $\Vert x\Vert _{_{E}} = \Vert y\Vert _{_{E}}$.
\end{enumerate}
\noindent Here, $m$ denotes Lebesgue measure and  $x^*$ denotes
the non-increasing, right-continuous rearrangement of $x$ given by
$$
x^{*}(t)=\inf \{~s\ge 0:m (\{u\in [0,1]:\,\mid x(u)\mid >s\})\le
t~\},\quad t>0.
$$
For basic properties of r.i. spaces, we refer to the monographs
\cite{KPS,LT-II}.  We note that for any r.i. space $E$ we have:
$L_\infty [0,1]\subseteq E\subseteq L_1[0,1].$  We will also work
with a r.i. space $E(\Omega,{\mathcal {P}})$ of measurable
functions on a probability space $(\Omega,{\mathcal {P}})$ given
by
$$
E(\Omega,{\mathcal {P}}):=\{f\in L_1(\Omega,{\mathcal {P}}):f^*\in
E\}, \quad \|f\|_{E(\Omega,{\mathcal {P}})}:=\|f^*\|_X.
$$
Here, the decreasing rearrangement $f^*$ is calculated with
respect to the measure ${\mathcal {P}}$ on $\Omega$.

Recall that for $0<\tau<\infty$, the dilation operator
$\sigma_\tau$ is defined by setting
$$\sigma_\tau x(t)=\begin{cases}
x(t/\tau),\;0\leq t\leq\min(1,\tau) \\
0,\; \min(1,\tau)<t\leq 1.
\end{cases}
$$  The dilation operators $\sigma_\tau$ are bounded in every r.i.  space
$E$. Denoting the space of all linear bounded operators on a
Banach space $E$ by ${\mathcal L}(E)$, we set
\[
 \alpha_E:=\lim\limits_{\tau\to
 0}\frac{\ln\|\sigma_\tau\|_{{\mathcal L}(E)}}{\ln\tau},\quad
 \beta_E:=\lim\limits_{\tau\to
 \infty}\frac{\ln\|\sigma_\tau\|_{{\mathcal L}(E)}}{\ln\tau}.
\]
The numbers $\alpha_E$ and $\beta_E$ belong to the closed interval
$[0,1]$ and are called the Boyd indices of $E$.

The K\"othe dual $E^\times $ of an r.i. space $E$ on $[0,1]$
consists of all measurable functions $y$ for which
$$
\Vert y\Vert _{_{E^{\times }}}:= \sup \Big\{\int _0^1\vert
x(t)y(t)\vert\,dt:\ x\in E,\ \Vert x \Vert _{_{E}}\leq 1\,\Big\}
<\infty.
$$
If $E^*$ denotes the Banach dual of $E$, then $E^\times \subset
E^{*}$ and $E^\times =E^{*}$ if and only if $E$ is separable. An
r.i. space $E$ is said to have the {\it Fatou property} if
whenever $\{f_n\}_{n=1}^\infty\subseteq E$ and $f$ measurable on
$[0,1]$ satisfy $f_n\to f$ a.e. on $[0,1]$ and $\sup _n\Vert
f_n\Vert _{_{E }} <\infty $, it follows  that $f\in E$ and $\Vert
f\Vert _{_{E}}\leq \liminf _{n\to \infty }\Vert f_n\Vert _{_{E}}$.
It is well-known that an r.i. space $E$ has the Fatou property if
and only if the natural embedding of $E$ into its K\"othe bidual
$E^{\times\times}$ is a surjective isometry.

\bigskip
Let us recall some classical examples of r.i. spaces on $[0,1]$.
Denote by $\Psi$ the set of all increasing continuous concave
functions on $[0,1]$ with
$\varphi(0)
=0$. Each function $\varphi\in\Psi$
generates the Lorentz space $\Lambda(\varphi)$ (see e.g.
\cite{KPS}) endowed with the norm
\[\|x\|_{\Lambda(\varphi)}=\int\limits_0^1 x^*(t)d\varphi(t)\]
and the Marcinkiewicz space $M(\varphi)$ endowed with the norm
\[
\|x\|_{M(\varphi)}=\sup\limits_{0<\tau\leq
1}\frac{1}{\varphi(\tau)}\int\limits_0^\tau x^*(t)dt.
\]
The space $M(\varphi)$ is not separable, but the space
\[\left \{x\in M(\varphi):\:\lim\limits_{\tau\to 0}\frac{1}{\varphi(\tau)}
\int\limits_0^\tau x^*(t)dt=0\right \}\] endowed with the norm
$\|\cdot \|_{M(\varphi)}$ is a separable r.i. space (denoted
further as $(M(\varphi)_0$), which coincides with the closure of
$L_\infty$ in $(M(\varphi),\|\cdot \|_{M(\varphi)})$.

It is well known (see e.g. \cite[Section II.1]{KPS}) that
$$
\beta_{M(\varphi)}=1\Longleftrightarrow
\alpha_{\Lambda(\varphi)}=0\Longleftrightarrow\forall t\in
(0,1)\exists (s_n)_{n\ge1} \subseteq (0,1)\ :\
\lim_{n\to\infty}\frac{\varphi(ts_n)}{\varphi(s_n)}=1;
$$
$$
\alpha_{M(\varphi)}=0\Longleftrightarrow
\beta_{\Lambda(\varphi)}=1\Longleftrightarrow\forall
\tau\ge1\exists (s_n)_{n\ge1} \subseteq (0,1)\ :\
\lim_{n\to\infty}\frac{\varphi(s_n\tau)}{\varphi(s_n)}=\tau.
$$


If $M(t)$ is a convex increasing function on $[0,\infty)$ such
that $M(0)=0$, then
the Orlicz space $L_M$ on $[0,1]$ (see e.g. \cite{KPS, LT-II}) is
a r.i. space of all $x\in L_1[0,1]$ such that
\[\|x\|_{L_M}:=\inf\{\lambda :\lambda >0,\;\int\limits_{0}^{1}
M(|x(t)|/\lambda)dt\leq 1\}<\infty.\] The function
$N_p(u)=e^{u^p}-1$ is convex for $p\geq1$ and is equivalent to a
convex function for $0<p<1$ (see e.g. \cite{Br1994, AsSu2005}).
The space $L_{N_p}$, $0<p<\infty$ is customarily denoted $\exp
(L_p)$.


\subsection{The Kruglov property in r.i.\ spaces}
Let $f$ be a random variable on $[0,1]$. By $\pi(f)$ we denote the
random variable $\sum_{i=1}^N f_i$, where $f_i$'s are independent
copies of $f$ and $N$ is a Poisson random variable with parameter
$1$ independent of the sequence $\{f_i\}$.

{\bf Definition.}\quad {\sl An r.i. space $E$ is said to have the
Kruglov property, if and only if $f\in E\Longleftrightarrow
\pi(f)\in E.$}

This property has been studied by M. Sh. Braverman \cite{Br1994}
which uses some earlier probabilistic constructions of V.M.
Kruglov \cite{K} and in \cite{AsSu2005,AsSu2006-1, AsSu2006-2} via
an operator approach. It was proved in \cite{AsSu2006-2}, that an
r.i. space $E$ satisfies the Kruglov property if and only if for
every sequence of independent mean zero functions $\{f_n\}\in E$
the following inequality holds
\begin{equation}\label{independent to disjoint}
||\sum_{k=1}^nf_k||_E\leq const\cdot
||\sum_{k=1}^n\overline{f}_k||_{Z^2_E}.
\end{equation}
Here, $Z^2_E$ is an r.i. space on $(0,\infty),$ equipped with a
norm
$$||x||=||x^*\chi_{[0,1]}||_E+||x^*\chi_{[1,\infty)}||_{L_2}$$
and the sequence $\{\bar f_k\}_{k=1}^n\subseteq Z^2_{X}$ is a
sequence of disjoint translates of $\{f_k\}_{k=1}^n\subseteq X,$
that is, $\bar f_k(\cdot)=f_k(\cdot-k+1)$. Note that inequality
\eqref{independent to disjoint} has been proved earlier in
~\cite{JoSch1989} (see inequality~(3) there) under the more
restrictive assumption that $E\supseteq L_p$ for some $p<\infty$.
Clearly, the latter assumption holds if $\alpha_E>0$.

\section{Operators $A_n$, $n\ge 0$}
Let $\Omega$ be the segment $[0,1],$ equipped with the Lebesgue
measure. 
Let $E$ be an arbitrary rearrangement invariant space on $\Omega.$

For every $n\geq 1$, we consider the operator
$A_n:E(\Omega)\rightarrow
E(\underbrace{\Omega\times\Omega\times\cdots\times\Omega}_{2n\
times})$ given by
\begin{multline*}
A_nf=(f\otimes r)\otimes(1\otimes1)\otimes\cdots\otimes(1\otimes 1)+(1\otimes1)\otimes(f\otimes r)\otimes\cdots\otimes(1\otimes 1)+\cdots\\
\cdots +(1\otimes1)\otimes \cdots \otimes(1\otimes1)\otimes(f\otimes r),
\end{multline*}
where $r$ is centered Bernoulli random variable. For brevity, we
will also use the following notation
$$A_nf=(f\otimes r)_1+(f\otimes r)_2+\cdots+(f\otimes r)_n.$$
We set $A_0=0.$


The following theorem is the main result of the present section.

\begin{thm}\label{alternative}
 The following alternative is valid in an arbitrary r.i.\ space $E.$
\begin{enumerate}
\item[(i).] $||A_n||_{{\mathcal L}(E)}=n$ for every natural $n;$

\item[(ii).] There exists a constant $\frac12\leq q<1,$ such that
$||A_n||_{{\mathcal L}(E)}\leq const\cdot n^q$ for all
$n\in\mathbb{N}.$
\end{enumerate}
\end{thm}

\begin{proof}

Since for all $n,m\geq0,$ we have
\begin{equation}\label{additive}
||A_{n+m}||_{{\mathcal L}(E)}\leq ||A_n||_{{\mathcal
L}(E)}+||A_m||_{{\mathcal L}(E)},
\end{equation}
and since $||f\otimes r||_E=||f||_E,$ we infer that
$||A_n||_{{\mathcal L}(E)}\leq n.$

Observing that $A_{mn}(f)$ and $A_m(A_n(f))$ are identically
distributed, we have
$$||A_{mn}(f)||_E=||A_m(A_n(f))||_E,\quad f\in E(\Omega).$$
Here, we identify  the element $A_nf\in
E(\Omega\times\cdots\times\Omega)$ with an element from
$E(\Omega)$ via a measure preserving transformation
$\underbrace{\Omega\times\cdots\times\Omega}_{2n\
times}\rightarrow\Omega.$ Hence,
\begin{equation}\label{multiplicative}||A_{mn}||_{{\mathcal L}(E)}\leq ||A_m||_{{\mathcal L}(E)}
\cdot||A_n||_{{\mathcal L}(E)}.\end{equation} Thus, we have the
following alternative:

\begin{enumerate}
\item[(i).] $||A_n||_{{\mathcal L}(E)}=n$ for every natural $n;$

\item[(ii).] There exists $n_0\geq2,$ such that
$||A_{n_0}||_{{\mathcal L}(E)}<n_0.$
\end{enumerate}

To finish the proof of Theorem \ref{alternative}, we need only to
consider the second case. Suppose there exists a constant
$\frac12\leq q<1,$ such that $||A_{n_0}||_{{\mathcal L}(E)}\leq
n_0^q.$ By \eqref{multiplicative} we have
$$||A_{n_0^m}||_{{\mathcal L}(E)}\leq ||A_{n_0}||_{{\mathcal L}(E)}^m\leq n_0^{qm},\ \forall m\in\mathbb{N}.$$

Every $n$ can be written as $\sum_{i=0}^ka_in_0^i,$ where $0\leq
a_i\leq n_0-1$ and $a_k\neq0.$ So, using \eqref{additive}  and
\eqref{multiplicative}, we have
$$||A_n||_{{\mathcal L}(E)}\leq \sum_{i=0}^k||A_{a_in_0^i}||_{\mathcal{L}(E)}\leq
\sum_{i=0}^k||A_{a_i}||_{\mathcal{L}(E)}n_0^{qi}\leq$$
$$\leq (\sum_{i=0}^k n_0^{qi})\max\limits_{1\leq s\leq n_0}\{||A_s||_{{\mathcal L}(E)}\}
\leq \frac{n_0^q\cdot n_0^{qk}}{n_0^q-1}\max\limits_{1\leq s\leq
n_0}\{||A_s||_{{\mathcal L}(E)}\}.$$ Now, using the fact that
$q\geq\frac12$ and $n_0\geq2,$ we have $n_0^q-1\geq(\sqrt{2}-1).$
So,
$$\frac1{n_0^q-1}\leq \sqrt{2}+1.$$
Since $n_0^k\leq n,$ we have
$$||A_n||_{{\mathcal L}(E)}\leq (\sqrt{2}+1)\cdot n_0^q
\cdot\max\limits_{1\leq s\leq n_0}\{||A_s||_{{\mathcal
L}(E)}\}\cdot n_0^{qk}\leq const\cdot n^q.$$ This proves the
theorem. \end{proof}

\begin{rem}\label{connection}
We record here an important connection between the estimates given
in Theorem \ref{alternative}(ii) above and the set $\Gamma_{\rm
iid}(E)$, where the r.i. space $E$ is separable. For $\frac12\leq
q\leq 1$ the following conditions are equivalent
\begin{itemize}
\item[(i)] $ ||A_n||_{{\mathcal L}(E)}\leq const\cdot n^q$, $n\ge
1$; \item[(ii)]\quad $ \frac{1}{q} \in \Gamma_{\rm iid}(E)$.
\end{itemize}
Indeed, the implication $(i)\Rightarrow(ii)$ is obvious. Now, let
the probability space $(\Omega, \mathcal{P})$ be the infinite
direct product of measure spaces $([0,1],m)$. Fix $f\in E$ and
consider the sequence $\{(f\otimes r)_n\}_{n\ge 1}\subset
E(\Omega, \mathcal{P})$. It follows from \cite[Lemma 3.4]{SeSu}
that this sequence is weakly null in $E(\Omega, \mathcal{P})$.
Since the spaces $E$ and $E(\Omega, \mathcal{P})$ are isometric,
we obtain the implication $(ii)\Rightarrow(i)$ via an application
of the uniform boundedness principle.
\end{rem}

We complete this section with an estimate of $\|A_n\|_{{\mathcal
L}(E)}$, $n\ge 1$ in general r.i. spaces with the Kruglov
property.

\begin{thm}\label{Kruglov} Let $E$ be a separable r.i. space.
If $\beta_E<1$ and if $E$ satisfies the Kruglov property, then
$||A_n||_{{\mathcal L}(E)}\leq const\cdot n^q$ for all
sufficiently large $n\ge 1$ and any $\beta_E<q<1$.
\end{thm}

\begin{proof}

It is proved in \cite[Proposition 2.2]{AsSeSu2007} (see also
\cite[Theorem~1]{MS2002}), that for every r.i. space $E$ and  an
arbitrary sequence of independent random variables
$\{f_k\}_{k=1}^n$ $(n\ge 1)$ from $E$, the right hand side of
\eqref{independent to disjoint} can be estimated as
\begin{equation}\label{quad sum}
||\sum_{k=1}^n\overline{f}_k||_{Z^2_E}\leq
6||(\sum_{k=1}^nf_k^2)^{\frac12}||_E.
\end{equation}

Now, assume in addition that the sequence $\{f_k\}_{k=1}^n$ $(n\ge
1)$ consists of independent identically distributed random
variables, $\|f_1\|_E=1$. Since $\beta_E<1,$ there exist $N$ and
$\beta_E<q<1$ such that for every $k\geq N$
$||\sigma_k||_{\mathcal{L}(E)}\leq k^{q}$. Fix $\varepsilon>0$
such that $\frac12+\varepsilon<q$. By \cite[Theorem 9]{SeSu}, in
every separable r.i. space $E$, the right hand side of \eqref{quad
sum} can be estimated as
\begin{equation}\label{quad bound}
||(\sum_{k=1}^nf_k^2)^{\frac12}||_E\leq\frac{4}{\varepsilon}\max_{1\leq
k\leq
n}(\frac{n}k)^{\frac12+\varepsilon}||\sigma_k||_{\mathcal{L}(E)}:=A,\quad
n\ge 1.
\end{equation}
So, the right hand side of \eqref{quad bound} can be estimated as

\begin{multline}\label{another_bound}
A\leq\frac{4}{\varepsilon}n^{\frac12+\varepsilon}
\max\{\max_{N\leq k\leq n}k^{-\frac12-\varepsilon}k^q,
\max_{1\leq k\leq N}k^{-\frac12-\varepsilon}||\sigma_k||_{\mathcal{L}(E)}\}=\\
=\frac{4}{\varepsilon}n^{\frac12+\varepsilon}
\max\{n^{q-\frac12-\varepsilon},const\}\leq const\cdot n^q.
\end{multline}
Recalling the definition of the operator $A_n$ and combining it
with \eqref{independent to disjoint}, \eqref{quad sum},
\eqref{quad bound}, \eqref{another_bound} yields the
assertion.\end{proof}

\begin{rem}
\begin{itemize}
\item[(i)] The assumption  $\beta_E<1$ in Theorem \ref{Kruglov} is necessary (see \cite[Theorem 4.2]{AsSeSu2005}).
For example, the space $E=L_1$ satisfies the Kruglov property and $\beta_E=1$. However,
$\|A_n\|_{\mathcal{L}(L_1)}= n$.
\item[(ii)] On the other hand, the condition that $E$ satisfies the Kruglov property is not optimal. In the following section, we will show that there are
Lorentz spaces which do not possess the Kruglov property and which
still satisfy the condition of Theorem \ref{alternative}(ii).
\end{itemize}
\end{rem}

\section{Operators $A_n$, $n\ge 1$ in Lorentz spaces.}

We need the following technical facts. The first lemma is elementary and its proof is omitted.
\begin{lem}\label{linearity} Let $\psi$ is a concave function on $[0,1].$ If there are points
$0\leq x_1\leq x_2\leq\cdots\leq x_n\leq1,$ such that
$$\frac1n(\psi(x_1)+\cdots \psi(x_n))=\psi(\frac1n (x_1+\cdots +x_n)),$$
then $\psi$ is linear on $[x_1,x_n].$
\end{lem}

\begin{lem}\label{estimate of expectation} Let $x_1,\cdots,x_n$ are independent random variables. The following inequality holds.
$$\mathbb{E}(|x_1+\cdots+x_n|)\leq \mathbb{E}(|x_1|)+\cdots+\mathbb{E}(|x_n|).$$
Moreover, the equality holds if and only if all $x_i'$s are simultaneously non-negative
(or non-positive).
\end{lem}

\begin{proof} We have
$$\mathbb{E}(|x_1|)+\cdots+\mathbb{E}(|x_n|)-\mathbb{E}(|x_1+\cdots+x_n|)=\mathbb{E}(|x_1|+\cdots+|x_n|-|x_1+\cdots+x_n|)\geq0.$$
By the independence of $x_i's,$ $i=1,2,\cdots,n$ we have
$sign(x_i),$ $i=1,2,\cdots,n$ are independent random variables. If
there exists a function $x_i,$ which is neither non-negative, nor
non-positive, then, for every other function $x_j,$ we have {\bf z-z-z-z}
\begin{multline*}
m(x_ix_j<0)=m(sign(x_i)>0,sign(x_j)<0)
+m(sign(x_i)<0,sign(x_j)>0)\\=m(sign(x_i)>0)m(sign(x_j)<0)
+m(sign(x_i)<0)m(sign(x_j)>0)>0.
\end{multline*}
Hence, there exists a set $A$ of positive measure such that $x_ix_j<0$ almost everywhere on $A.$ This guarrantees that $|x_1|+\cdots+|x_n|>|x_1+\cdots+x_n|$ almost everywhere on $A.$ This is sufficient for the strict inequality to hold.\end{proof}

We need to consider the following properties of the function $\psi.$
\begin{equation}\label{first property}a_\psi:=\limsup_{u\rightarrow0}\frac{\psi(ku)}{\psi(u)}<k.\end{equation}
\begin{equation}\label{second property}c_\psi:=\limsup_{u\rightarrow0}\frac{\psi(u^l)}{\psi(u)}<1.\end{equation}
\begin{equation}\label{general limit property}\limsup_{u\rightarrow0}\frac1{\psi(u)}\sum_{s=1}^n\psi(2^{1-s}\binom{n}{s}u^s)<n.\end{equation}

\begin{prop}\label{two limits into general} Suppose, there exist $k\geq2$ such that
\eqref{first property} holds and $l\geq2$ such that
\eqref{second property} holds. Then, \eqref{general limit
property} holds for all sufficiently large $n.$
\end{prop}

\begin{proof}
Consider the sum $\sum_{s=1}^n\psi(\binom{n}{s}2^{1-s}u^s).$ For any sufficiently large $n,$ we write
$$\sum_{s=1}^n=\sum_{s=1}^{1+[\frac{n}k]}+\sum_{s=2+[\frac{n}k]}^n.$$
Consequently, the upper limit in \eqref{general limit property}
can be estimated as
\begin{equation}\label{decomposition}\begin{split}
\limsup_{u\rightarrow0}\frac1{\psi(u)}\sum_{s=1}^n\psi(\binom{n}{s}2^{1-s}u^s)\leq\limsup_{u\rightarrow0}\frac1{\psi(u)}\sum_{s=1}^{1+[\frac{n}k]}\psi(\binom{n}{s}2^{1-s}u^s)+\\
+\limsup_{u\rightarrow0}\frac1{\psi(u)}\sum_{s=2+[\frac{n}k]}^n\psi(\binom{n}{s}2^{1-s}u^s)
\end{split}\end{equation}

Consider the first upper limit in \eqref{decomposition}. Since
$\psi$ is concave, we have
$$\sum_{s=1}^{1+[\frac{n}k]}\psi(\binom{n}{s}2^{1-s}u^s)\leq(1+[\frac{n}{k}])\psi(\frac1{1+[\frac{n}k]}\sum_{s=1}^{1+[\frac{n}k]}\binom{n}{s}2^{1-s}u^s)=$$
$$=(1+[\frac{n}{k}])\psi(\frac1{1+[\frac{n}k]}(nu+o(u)))\leq(1+[\frac{n}{k}])\psi(ku(1+o(1))).$$
Hence, the first upper limit in \eqref{decomposition} is bounded
from above by
$$(1+[\frac{n}k])a_{\psi}=n\cdot\frac{a_{\psi}}{k}+o(n).$$

Consider the second upper limit in \eqref{decomposition}. It is
clear that for all $\frac1k n\leq s\leq n$
$$\binom{n}{s}\cdot2^{1-s}\leq 2^n$$
and
$$\binom{n}{s}2^{1-s}u^s\leq 2^nu^{\frac1k n}=(2^ku)^{\frac1k n}.$$
Thus, the second upper limit in \eqref{decomposition} can be
estimated as
$$\limsup_{u\rightarrow0}\frac1{\psi(u)}\sum_{s=2+[\frac{n}k]}^n\psi(\binom{n}{s}2^{1-s}u^s)\leq n(1-\frac1k)\limsup_{u\rightarrow0}\frac{\psi((2^ku)^{\frac{n}{k}})}{\psi(u)}.$$
Substituting variable $w=2^ku$ on the right hand side, we have
$$n(1-\frac1k)\limsup_{w\rightarrow0}\frac{\psi(w^{\frac{n}{k}})}{\psi(2^{-k}w)}.$$
By the concavity of $\psi,$ we have
$\psi(2^{-k}w)\geq2^{-k}\psi(w).$ Therefore, the second upper
limit in \eqref{decomposition} is bounded from above by
$$n(1-\frac1k)2^k\limsup_{w\rightarrow0}\frac{\psi(w^{\frac{n}{k}})}{\psi(w)}.$$
Now, we observe that
\begin{equation}
\limsup_{w\rightarrow0}\frac{\psi(w^m)}{\psi(w)}\leq c_{\psi}^{\frac{\log(m)}{\log(l)}-1}.
\end{equation}
Indeed, let $l^r\leq m\leq l^{r+1},$
$$\frac{\psi(w^m)}{\psi(w)}\leq\frac{\psi(w^{l^r})}{\psi(w)}=\frac{\psi(w^{l^r})}{\psi(w^{l^{r-1}})}\cdots\frac{\psi(w^l)}{\psi(w)}$$
and
$$\limsup_{w\rightarrow0}\frac{\psi(w^m)}{\psi(w)}\leq c_{\psi}^r \leq c_{\psi}^{\frac{\log(m)}{\log(l)}-1}.$$
If $n$ tends to infinity, then, thanks to the assumption $c_{\psi}<1,$ we have $$n(1-\frac1k)2^k\limsup_{w\rightarrow0}\frac{\psi(w^{\frac{n}k})}{\psi(w)}=o(n).$$

Therefore, the upper limit in \eqref{general limit property} (see
also \eqref{decomposition}) is bounded from above by
$$\frac{a_{\psi}}k n+o(n).$$
Thus, the upper limit in \eqref{general limit property} is
strictly less then $n$ for every sufficiently large $n.$
\end{proof}

Let the function $g_n$ be defined by
\begin{equation}\label{g definition}
g_n(u):=\frac{||A_n\chi_{[0,u]}||_{\Lambda(\psi)}}{n||\chi_{[0,u]}||_{\Lambda(\psi)}}
=\frac1{n\psi(u)}\sum_{s=1}^n\psi(m(|(\chi_{[0,u]}\otimes
r)_1+\cdots+(\chi_{[0,u]}\otimes r)_n|\geq s)).
\end{equation}
It is obvious that $0\leq g_n\leq 1$.

\begin{rem} The second equality in \eqref{g definition} is a corollary of
\cite[II.5.4]{KPS}.
\end{rem}

\begin{prop}\label{tech fignya} For sufficiently large $n,$ we have $g_n(u)<1$ for all
$u\in(0,1]$.
\end{prop}

\begin{proof}  Since $\psi$ is concave, we have
\begin{multline}\label{first bound of g}
\sum_{s=1}^n\psi(m(|(\chi_{[0,u]}\otimes r)_1+\cdots+(\chi_{[0,u]}\otimes r)_n|\geq s))\leq \\
\leq n\cdot\psi(\frac1n\sum_{s=1}^nm(|(\chi_{[0,u]}\otimes
r)_1+\cdots+(\chi_{[0,u]}\otimes r)_n|\geq s)).
\end{multline}
Note, that if random variable $\xi_n$ takes the values $0,1,\cdots,n$ then
\begin{equation}\label{expectation appears}
\sum_{s=1}^nm(\xi_n\geq s)=\mathbb{E}(\xi_n).
\end{equation}
By \eqref{expectation appears}, the right-hand side of
\eqref{first bound of g} is equal to $n\psi(\frac1n
\mathbb{E}(|(\chi_{[0,u]}\otimes r)_1+\cdots+(\chi_{[0,u]}\otimes
r)_n|)).$ By Lemma \ref{estimate of expectation}, we have
\begin{equation}\label{strict estimate of expectation}
\frac1n \mathbb{E}(|(\chi_{[0,u]}\otimes
r)_1+\cdots+(\chi_{[0,u]}\otimes
r)_n|)<\mathbb{E}(|\chi_{[0,u]}\otimes r|)=u.
\end{equation}
Taking $\psi$, we obtain
\begin{equation}\label{second bound of g}
n\psi(\frac1n\mathbb{E}(|(\chi_{[0,u]}\otimes
r)_1+\cdots+(\chi_{[0,u]}\otimes r)_n|))\leq
n\psi(\mathbb{E}(|\chi_{[0,u]}\otimes r|)).
\end{equation}
The right hand side of \eqref{second bound of g} is equal to
$n\psi(u).$

Let us assume that $g_n(u)=1,$ for some $u>0$ and some $n>1.$ It
then follows, that both inequalities \eqref{first bound of g} and
\eqref{second bound of g} are actually equalities.

The equality
\begin{multline*}
\sum_{s=1}^n\psi(m(|(\chi_{[0,u]}\otimes r)_1+\cdots+(\chi_{[0,u]}\otimes r)_n|\geq s))=\\
=n\cdot\psi(\frac1n\sum_{s=1}^nm(|(\chi_{[0,u]}\otimes
r)_1+\cdots+(\chi_{[0,u]}\otimes r)_n|\geq s))
\end{multline*}
implies, by Lemma \ref{linearity}, that $\psi$ is linear on the
interval $[a_1,b_1]$ with $a_1=m(|(\chi_{[0,u]}\otimes
r)_1+\cdots+(\chi_{[0,u]}\otimes r)_n|\geq n),$ and
$b_1=m(|(\chi_{[0,u]}\otimes r)_1+\cdots+(\chi_{[0,u]}\otimes
r)_n|\geq1).$

Since the inequality in \eqref{second bound of g} is actually an
equality, we derive from \eqref{strict estimate of expectation}
and \eqref{second bound of g}, that  $\psi$ must be a constant on
the interval $[a_2,b_2]$ with $a_2=\frac1n
\mathbb{E}(|(\chi_{[0,u]}\otimes r)_1+\cdots+(\chi_{[0,u]}\otimes
r)_n|),$ and $b_2=\mathbb{E}(|\chi_{[0,u]}\otimes r|)].$ Since
$\psi$ is increasing and concave function, it must be a constant
on $[a_2,1].$

Since, by \eqref{expectation appears},
$$\frac1n \mathbb{E}(|(\chi_{[0,u]}\otimes r)_1+\cdots+(\chi_{[0,u]}\otimes r)_n|)>m(|(\chi_{[0,u]}\otimes r)_1+\cdots+(\chi_{[0,u]}\otimes r)_n|\geq n)$$
and
$$\frac1n \mathbb{E}(|(\chi_{[0,u]}\otimes r)_1+\cdots+(\chi_{[0,u]}\otimes r)_n|)<m(|(\chi_{[0,u]}\otimes r)_1+\cdots+(\chi_{[0,u]}\otimes r)_n|\geq1)),$$
we have $a_1<a_2<b_2.$ So, the intersection of the intervals $[a_1,b_1]$ and $[a_2,1]$ contains an interval $[a_3,b_3]$ with $a_3<b_3.$

Since $\psi$ is a linear function on the $[a_1,b_1]$ and is a constant on the $[a_2,1]$ it must be a constant on $[a_1,1]$ that is on the interval
$$[m(|(\chi_{[0,u]}\otimes r)_1+\cdots+(\chi_{[0,u]}\otimes r)_n|\geq n),1]=[2^{1-n}u^n,1].$$
Thus, $\psi$ is a constant on the interval $[2^{1-n},1]\subset [2^{1-n}u^n,1],$ which is not the case for sufficiently large $n.$ So, $g_n(u)<1$ for all sufficiently large $n.$
\end{proof}

\begin{lem}\label{limit equivalence} For the function $g_n,$ defined in Proposition \ref{tech fignya}, we have
$$\limsup_{u\rightarrow0}g_n(u)=\limsup_{u\rightarrow0}\frac1{n\psi(u)}\sum_{s=1}^n\psi(2^{1-s}\binom{n}{s}u^s).$$
\end{lem}
\begin{proof}For every $s\geq 1,$ using a formula for conditional probabilities, we have
\begin{equation*}
m(|(\chi_{[0,u]}\otimes r)_1+\cdots+(\chi_{[0,u]}\otimes r)_n|\geq
s)=\sum_{k=1}^n\binom{n}{s}u^k(1-u)^{n-k}m(|r_1+\cdots+r_k|\geq s).
\end{equation*}
Actually, the summation above is taken from $k=s$ up to $n,$ since $m(|r_1+\cdots+r_k|\geq s)=0$ for every $k<s.$\\
If now $u\rightarrow0,$ then, for every $s\geq1$ and $k>s,$ we have $\binom{n}{k}u^k(1-u)^{n-k}=o(u^s).$ Therefore,
\begin{equation}\label{main term selection}
m(|(\chi_{[0,u]}\otimes r)_1+\cdots+(\chi_{[0,u]}\otimes r)_n|\geq
s)=2^{1-s}\binom{n}{s}u^s(1+o(1)).
\end{equation}
Since $\psi$ is concave, we have
\begin{equation}\label{sublinearity}
\psi(\frac1{m}u)\leq\frac1{m}\psi(u),\ 0<m\leq1.
\end{equation}
This implies
\begin{equation}\label{psi limit property}
\lim\limits_{u\rightarrow0}\frac{\psi(u(1+o(1)))}{\psi(u)}=1.
\end{equation}
After applying \eqref{main term selection} and \eqref{psi limit
property} to the definition of $g_n$ in \eqref{g definition}, we
obtain the assertion of the lemma.
\end{proof}

The following theorem is the main result in this section.
\begin{thm}\label{characterization} Let $\psi\in\Psi.$ The following conditions are equivalent.

$(i)$ $||A_n||_{\mathcal{L}(\Lambda(\psi))}<n$ for all
sufficiently large $n;$

$(ii)$ Estimates \eqref{first property} and \eqref{second
property} hold for some $k\geq2$ and $l\geq2$.
\end{thm}

\begin{rem} Note that condition $(i)$ above is equivalent to the assumption that
$||A_{n_0}||_{\mathcal{L}(\Lambda(\psi))}<n_0$ for some $n_0>1$
(see Theorem \ref{alternative}).
\end{rem}

\begin{proof} We are interested whether there exist $n\in\mathbb{N}$ and $c<n,$ such that
\begin{equation}\label{an bound}
||A_nf||_{\Lambda(\psi)}\leq c||f||_{\Lambda(\psi)},\
\ f\in\Lambda(\psi).
\end{equation}

We will use the following known description of extreme points of
the unit ball in $\Lambda(\psi).$ A function $f\in
extr(B_{\Lambda(\psi)}(0,1))$ if and only if
$$|f|=\frac{\chi_{A}}{||\chi_{A}||_{\Lambda(\psi)}}$$
for some measurable set $A\subset [0,1].$ Here $\chi_A$ is the
indicator function of the set $A$. This means that $f$ is of the
form
$$f=\frac{\chi_{A_1}-\chi_{A_2}}{\psi(m(A_1\cup A_2))}$$
with $A_1$ and $A_2$ having empty intersection. It is sufficient
to verify \eqref{an bound} only for functions $f$ as above (see
\cite[Lemma II.5.2]{KPS}).

Clearly, $f\otimes r$ and $|f|\otimes r$ are identically
distributed random variables. Therefore, $A_n(f)$ and $A_n(|f|)$
are also identically distributed ones. Furthermore,
$||A_m(f)||=||A_m(|f|)||$ and $||f||=||\,|f|\,||.$ Thus, we need
to check \eqref{an bound} for indicator functions only. It is
sufficient to take $A$ of the form $[0,u],$ $0<u\leq1.$

Using the notation $g_n(\cdot)$ introduced in \eqref{g
definition}, we see that \eqref{an bound} is equivalent to
\begin{equation}\label{gn bound}
\sup_u g_n(u)<1.
\end{equation}

Now, we are ready to finish the proof of the theorem.

[Necessity] Fix $n$ such that
$||A_n||_{\mathcal{L}(\Lambda(\psi))}<n.$ It follows from the
argument above that \eqref{gn bound} holds. Now, we immediately
infer from Lemma \ref{limit equivalence} and the definition of
$g_n(\cdot)$ that
$$\limsup_{u\rightarrow0}\frac1{n\psi(u)}\sum_{s=1}^n\psi(\binom{n}{s}2^{1-s}u^s)<1,$$
which is equivalent to
 \eqref{general limit property}. Thus,
$$\limsup_{u\rightarrow0}\frac{\psi(nu)}{n\psi(u)}
=\limsup_{u\rightarrow0}\frac{\psi(2^{1-1}\binom{n}{1}u^1)}{n\psi(u)}
\leq\limsup_{u\rightarrow0}\frac1{n\psi(u)}\sum_{s=1}^n\psi(2^{1-s}\binom{n}{s}u^s)<1.$$

Suppose that \eqref{second property} fails. We have
$$\limsup_{u\rightarrow 0}\frac{\psi(u^l)}{\psi(u)}=1$$
for every $l\geq 1.$ Since $\binom{n}{s}2^{1-s}u^s\geq u^{n+1}$ for every $s=1,2,\cdots,n$
and every sufficiently small $u,$ we have
$$\limsup_{u\rightarrow 0}\frac1{n\psi(u)}\sum_{s=1}^n\psi(\binom{n}{s}2^{1-s}u^s)
\geq\limsup_{u\rightarrow 0}\frac{n\psi(u^{n+1})}{n\psi(u)}=1.$$
This contradicts with \eqref{general limit property} and completes
the proof of necessity.

[Sufficiency] Fix $k\geq2$ (respectively, $l\geq2$) such that
\eqref{first property} (respectively, \eqref{second property})
holds. Then, for sufficiently large $n,$ \eqref{general limit
property} also holds. By Lemma \ref{limit equivalence}, we have
\begin{equation}\label{gn limit bound}
\limsup_{u\rightarrow0}g_n(u)<1
\end{equation}
for all sufficiently large $n.$ By Proposition \ref{tech fignya},
we have $g_n(u)<1$ for all sufficiently large $n$ and for all
$u\in(0,1].$ Therefore, by \eqref{gn limit bound}, \eqref{gn
bound} holds for sufficiently large $n.$ Then (see the argument at
the beginning of the proof),
$||A_n||_{\mathcal{L}(\Lambda(\psi))}<n$ for sufficiently large
$n.$
\end{proof}

Combining Theorems \ref{alternative} and \ref{characterization},
we have

\begin{cor}\label{lorentz alternative} For every function $\psi,$ one of
the following two mutually excluding
alternatives holds.

\begin{enumerate}
\item  There exist $q\in[\frac12,1)$ and $C>0,$ such that the
operator $A_n:\Lambda(\psi)\rightarrow\Lambda(\psi)$ satisfies
$$||A_n||_{\mathcal{L}(E)}\leq C\cdot n^q,\ n\geq1.$$

\item Either for every $k\in\mathbb{N},$
\begin{equation}\label{first bad condition}
\limsup_{u\rightarrow0}\frac{\psi(ku)}{\psi(u)}=k
\end{equation}
or for every $l\in\mathbb{N},$
\begin{equation}\label{second bad condition}
\limsup_{u\rightarrow0}\frac{\psi(u^l)}{\psi(u)}=1.
\end{equation}
\end{enumerate}
\end{cor}

\begin{rem} \begin{itemize}\item[(i)] The condition \eqref{first bad
condition} is equivalent to the assumption
$\beta_{\Lambda(\psi)}=1$. \item[(ii)] The condition \eqref{second
bad condition} implies (but not equivalent to) the condition
$\alpha_{\Lambda(\psi)}=0$. 
In the last
section of this paper, we will present an example $\psi\in \Psi$
failing \eqref{second bad condition}
 such that the Lorentz space $\Lambda(\psi)$  fails the Kruglov property.
\end{itemize}
\end{rem}

\section{Operators $A_n$, $n\ge 1$ in Orlicz spaces $exp(L_p)$.}

The space $exp(L_p)$ satisfies Kruglov property if and only if
$p\leq1$ (see \cite{Br1994,AsSu2005}). The space $exp(L_p)$ is
2-convex for all $0<p<\infty$ (see e.g. \cite[1.d]{LT-II}). Now,
we immediately infer from \cite{AsSeSu2007} that $\Gamma_{\rm
iid}(exp(L_p)_0)=\Gamma_{\rm i}(exp(L_p)_0)=[1,2]$ for all
$0<p\leq1$ (here, $exp(L_p)_0$ is the separable part of the space
$exp(L_p)$). Using Remark \ref{connection}, we have
$||A_n||_{\mathcal{L}(exp(L_p)_0)}\leq const\cdot n^{\frac12}$ for
all $n\geq1$ and $0<p\leq1$. It easily follows that in fact,
$||A_n||_{\mathcal{L}(exp(L_p))}\leq const\cdot n^{\frac12}$ for
all $n\geq1$ and $0<p\leq1$.  In this section, we prove the
estimate $||A_n||_{\mathcal{L}(exp(L_p))}\leq const\cdot
n^{\frac12}$ (respectively, $||A_n||_{\mathcal{L}(exp(L_p))}\leq
const\cdot n^{1-1/p}$) for all $1<p\leq 2$ (respectively, $2\leq
p<\infty$.) To this end, it is convenient to view $exp(L_p)$ as a
Marcinkiewicz space $M(\psi_p)$ with
$\psi_p(t)=t\log^{\frac1p}(\frac{e}t)$ (see \cite[Lemma
4.3]{AsSu2005}). The following simple lemma is crucial.

\begin{lem}\label{expl2 normal}
There exists $\Psi\ni\psi\sim\psi_2,$ such that the random
variable $\psi'\otimes r$ is Gaussian.
\end{lem}

\begin{proof} Setting $F(t):=\frac2{\sqrt{\pi}}\int_t^{\infty}e^{-z^2}dz$, $t\ge 0$
and denoting its inverse by $G$, we clearly have that $G\otimes r$
is Gaussian. From the obvious inequality
$$c_1\cdot e^{-2t^2}\leq F(t)\leq c_2 e^{-t^2},$$
substituting $t=G(z),$ we obtain
$$c_1\cdot e^{-2G^2(z)}\leq z\leq c_2 e^{-G^2(z)}$$
or, equivalently,
$$\frac1{\sqrt{2}}\log^{\frac12}(\frac{c_1}z)\leq G(z)\leq \log^{\frac12}(\frac{c_2}z).$$
This means
$$\psi(t)=\int_0^tG(z)dz\sim\int_0^t\log^{\frac12}(\frac{e}z)dz\sim t\log^{\frac12}(\frac{e}t)=\psi_2(t).$$
\end{proof}

\begin{thm}\label{explp bound}
\begin{itemize}
\item[(i)] For every $1\leq p\leq 2,$ we have
$||A_n||_{\mathcal{L}(exp(L_p))}\leq const\sqrt{n}.$

\item[(ii)]For every $2\leq p\leq\infty,$ we have
$||A_n||_{\mathcal{L}(exp(L_p))}\leq const\cdot n^{1-1/p}.$
\end{itemize}
\end{thm}
\begin{proof} (i).\quad By Lemma \ref{expl2 normal} $exp(L_2)=M(\psi)$, $\psi\in
\Psi$ where $\psi'\otimes r$ is Gaussian. Recall the following
description of the extreme points of the unit ball in
Marcinkiewicz spaces (see \cite{Ryff}): a function $f$ is an
extreme point of the unit ball in $M({\psi})$ if and only if
$f^*=\psi'$. Since $||A_nx||_{M({\psi})}=||A_n\psi'||_{M({\psi})}$
for any $x\in M_{\psi}$ with $x^*=\psi'$, we infer that
$||A_n\psi'||_{M({\psi})}=||A_n||_{\mathcal{L}(M({\psi}))}$, $n\ge
1$. Since the $\psi'\otimes r$ is Gaussian, the function
$$\frac{(\psi'\otimes r)_1+\cdots+(\psi'\otimes r)_n}{\sqrt{n}}$$
is also Guassian, in particular, its rearrangement coincides with
 $\psi'$. This means $||A_n||_{\mathcal{L}(M_{\psi})}=\sqrt{n}$.
The result now follows by interpolation between $exp(L_{1})$ and
$exp(L_{2})$, since for every $0<p_1\leq p_2\leq\infty$ we have
$$[exp(L_{p_1}),exp(L_{p_2})]_{\theta,\infty}=exp(L_p)$$
with $\frac1p=\frac{1-\theta}{p_1}+\frac{\theta}{p_2}$ (see, for
example \cite{BrKrug}).

(ii).\quad Noting that $||A_n||_{\mathcal{L}(L_\infty)}=n$, $n\ge
1$, the assertion follows from (i) by applying the real method of
interpolation to the couple $(exp(L_{2}),L_\infty)$ as above.
\end{proof}

\section{Applications to Banach-Saks index sets}

The first main result of this section  characterizing a subclass
of the class of all r.i. spaces $E$ such that $\Gamma_{\rm
iid}(E)= \Gamma_{\rm i}(E)$ is given in Theorem \ref{firstmain}
below. We firstly need a modification of the subsequence splitting
result from \cite[Theorem 3.2]{new-16-Sem-Suk}. We present
necessary details of the proof for convenience of the reader.

\begin{thm}\label{firsttechnical} Let $\{x_n\}_{n\ge 1}$ be a weakly null sequence of independent functions in a separable r.i. space $E$ with the Fatou property. Then, there exists a subsequence $\{y_n\}_{n\ge 1}\subset\{x_n\}_{n\ge 1},$ which can be split as $y_n=u_n+v_n+w_n, n\ge 1$. Here $\{u_n\}_{n\ge 1}$ is a weakly null sequence of independent identically distributed functions, the sequence $\{v_n\}_{n\ge 1}$ is also weakly  null and consists of the elements with pairwise disjoint support and $\|w_n\|_E\to 0$ as $n\to \infty$.
\end{thm}

\begin{proof} Let the probability space $(\Omega, \mathcal{P})$ be the infinite direct product of measure spaces $([0,1],m)$. Without loss of generality, we assume that $E=E(\Omega)$ and that each function $x_n$
depends only on the $n-$th coordinate. That is the
following holds
$$x_n=\underbrace{1\otimes\cdots\otimes1}_{(n-1)\
times}\otimes h_n\otimes 1\otimes\cdots,\quad h_n\in E(0,1),\ \quad n\ge 1.$$ Consider the sequence
$\{g_n\}_{n\ge 1}=\{h^*_n\}_{n\ge 1}\subset E(0,1)$. Since
$$||x_n||_E=||g_n||_E\geq ||g_n\chi_{[0,s]}||_E\geq g_n(s)||\chi_{[0,s]}||_E,\quad s\in [0,1]$$
and the sequence $\{x_n\}$ is bounded, it follows from Helly Selection theorem that there exists a
subsequence $\{g_n^1\}\subset\{g_n\},$ which converges almost
everywhere on $[\frac12,1]$. Repeating the argument, we get a
subsequence $\{g_n^2\}\subset\{g_n^1\},$ which converges almost
everywhere on $[\frac13,1],$  etc. Thus, there exists a function
$h\in L_1(0,1)$ to which the diagonal sequence $\{g_n^n\}_{n\ge
1}=\{(h_n^n)^*\}_{n\ge 1}$ converges almost everywhere. The Fatou
property of $E$ guarantees that $h\in E(0,1)$ and $\|h\|_E\leq 1$.
There is an operator $P_n:L_1(0,1)\to L_1(0,1)$ of the form
$(P_nx)(t)=\alpha(t)x(\gamma(t))$ (here $|\alpha(t)|=1$ and
$\gamma$ is a measure preserving transformation of the interval
$(0,1)$ into itself), such that $P_ng_n^n=h_n^n$, $n\ge 1$ (see
e.g. \cite{KPS}). Now, put
$$y_n:=1\otimes1\otimes\cdots\otimes1\otimes h^n_n\otimes1\cdots,\quad n\ge 1,$$
$$u_n:=1\otimes1\otimes\cdots\otimes1\otimes (P_n h)\otimes1\cdots,\quad n\ge 1.$$
It is clear, that functions $u_n$ are independent. {\bf z-z-z-z} The proof is
finished by repeating the remaining argument from \cite[Theorem
3.2]{new-16-Sem-Suk}.
\end{proof}

\begin{thm}\label{firstmain} For an arbitrary separable r.i. space $E$ with the Fatou property,
we have
$$
\Gamma_{\rm iid}(E)= \Gamma_{\rm i}(E)\Longleftrightarrow
\Gamma_{\rm iid}(E)\subseteq \Gamma_{\rm d}(E).
$$
\end{thm}

\begin{proof} If $\Gamma_{\rm iid}(E)= \Gamma_{\rm i}(E)$, then
the embedding $\Gamma_{\rm iid}(E)\subseteq \Gamma_{\rm d}(E)$
follows immediately from \cite [Lemma 4.1(ii)]{AsSeSu2005}.
Suppose now that $\Gamma_{\rm iid}(E)\subseteq \Gamma_{\rm d}(E)$
and let $\{f_k\}_{k\ge 1}\subset E$ be a normalized weakly null
sequence of independent random variables on $[0,1]$. Passing to a
subsequence and applying the preceding theorem, we may assume that
$f_n=u_n+v_n+w_n, n\ge 1$, where $\{u_n\}_{n\ge 1}$ is a weakly
null sequence of independent identically distributed functions,
the sequence $\{v_n\}_{n\ge 1}$ is also weakly  null and
 consists of the elements with pairwise disjoint support and $\|w_n\|_E\to 0$ as $n\to
 \infty$. Due to the latter convergence, we may assume without loss of generality that
 $||w_k||_{E}\leq 2^{-k}$ and so for every subsequence $\{z_n\}\subset\{w_n\},$ we have
$$||\sum_{k=1}^nz_k||_{E}\leq 1.$$
If, in addition,  $\frac {1}{q} \in \Gamma_{\rm iid}(E)$, then our
assumptions also guarantee that there are constants $C_2, C_3>0$
$$||\sum_{k=1}^nu_k||_{E}\leq C_2\cdot n^{q},\quad
||\sum_{k=1}^nv_k||_{E}\leq C_3\cdot n^{q}.$$
\end{proof}

We will illustrate the result above in the settings of:
$(\alpha)$ r.i. spaces satisfying an upper $2$-estimate; $(\beta)$
Lorentz spaces $\Lambda(\varphi)$ and Marcinkiewicz spaces
$M(\varphi)_0$, $\varphi\in \Psi$; and $(\gamma)$ classical
$L_{p,q}$-spaces.

$(\alpha)$\quad Recall that a Banach lattice $X$ is said to
satisfy an \textit{upper\ }$2-$\textit{estimate}, if there exists
a constant $C>0$ such that for every finite sequence
$(x_{i})_{_{i=1}}^{n}$ of pairwise disjoint elements in $X$
\begin{equation*}
\left\Vert  \sum_{j=1}^{n}x_{j}\right\Vert _{X}\leq C\left(
\sum_{j=1}^{n}\Vert x_{j}\Vert _{X}^{2}\right) ^{1/2}.
\end{equation*}%

\begin{cor}\label{upper2estimate} If $E$ is a separable r.i. space with the Fatou
property and satisfying an upper $2$-estimate, then $ \Gamma_{\rm
iid}(E)= \Gamma_{\rm i}(E)$.
\end{cor}

\begin{proof} The assumption that the space $E$ satisfies an upper
$2$-estimate implies immediately that $2\in \Gamma_{\rm d}(E)$ and
hence $[1,2]\subseteq \Gamma_{\rm d}(E)$. Noting that $
\Gamma_{\rm iid}(E) \subseteq [1,2]$ (see  \cite [Lemma
4.1(i)]{AsSeSu2005} the result now follows from Theorem
\ref{firstmain}.
\end{proof}

$(\beta)$\quad Although Lorentz spaces do not satisfy an upper
$2$-estimate, we have $$\Gamma_{d}(\Lambda(\psi))=[1,\infty)$$
(see e.g. the proof of \cite [Corollary 4.8]{AsSeSu2005}) and
similarly, $\Gamma_{d}(M(\psi)_0)=[1,\infty)$ (see e.g. \cite
[p.897]{AsSeSu2005}) for any $\psi\in \Psi$. Although, the
Marcinkiewicz spaces $(M(\psi)_0)$ do not possess the Fatou
property, applying the modification of Theorem
\ref{firsttechnical} similar to to \cite[Lemma 3.6]{AsSeSu2005},
we obtain the following corollary from Theorem \ref{firstmain}.

\begin{cor}\label{coincidence_in_Lor} For every $\psi\in \Psi$, we have $\Gamma_{\rm
i}(\Lambda(\psi))= \Gamma_{\rm iid}(\Lambda(\psi))$ and
$\Gamma_{\rm i}(M(\psi)_0)= \Gamma_{\rm iid}(M(\psi)_0)$.
\end{cor}


$(\gamma)$\quad We will now show that the equality $\Gamma_{\rm
i}(E)=\Gamma_{\rm iid}(E)$ fails in the important subclass of r.i.
space which plays a significant role in the interpolation theory
\cite{KPS,LT-II}. Recall the definition of the Lorentz spaces
$L_{p,q}$, $1<p,q<\infty$: $x\in L_{p,q}$ if and only if the
quasi-norm
\[
 \|x\|_{p,q}=
   \dfrac{q}{p}\left(\displaystyle\int\limits_0^1
   \left(x^*(t)t^{1/p}\right)^q\dfrac{dt}{t}\right)^{1/q}, 
\]
is finite. The expression $\|\cdot \|_{p,q}$ is a norm if
$1\leqslant q\leqslant p$ and is equivalent to a (Banach) norm if
$q>p$.

We will now show that $\Gamma_{\rm i}(L_{p,q})\neq\Gamma_{\rm
iid}(L_{p,q})$, provided $1<q<p<2$. To this end, we firstly
observe that every normalized sequence $\{v_n\}_{n\ge 1}\subset
L_{p,q}$ of functions with disjoint support contains a subsequence
spanning the space $l_q$ (see \cite[Lemma 2.1]{CD}). In
particular, $\Gamma_{\rm d}(L_{p,q})\subset\Gamma(l_q)=[1,q]$ and
so, by \cite[Lemma 4.1(ii)]{AsSeSu2005}, we have $\Gamma_{\rm
i}(L_{p,q})\subseteq [1,q]$. Next, it is proved in \cite[Corollary
5.2]{Br1996} (see also \cite{CarDil89}) {\bf z-z-z-z} that if $p<2$ then for every sequence of identically
distributed independent random variables we have
$$||\sum_{k=1}^nx_k||_{L_{p,q}}=o(n^{\frac1p}),$$
which implies, in particular, that $[1,p]\subseteq \Gamma_{\rm
iid}(L_{p,q})$. This shows that $(q,p]\subseteq \Gamma_{\rm
iid}(L_{p,q})\setminus \Gamma_{\rm i}(L_{p,q})$ as soon as
$1<q<p<2$.

Our second main result in this section completely characterizes
the subclass of all Lorentz spaces $\Lambda(\psi)$, $\psi\in \Psi$
whose Banach-Saks index set $\Gamma_{\rm i}(\Lambda(\psi))$ is
non-trivial.

\begin{thm}\label{mainsecond} $\Gamma_{\rm iid}(\Lambda(\psi))\neq\{1\}$ if and only if
the function $\psi$  satisfies conditions \eqref{first property}
and \eqref{second property} for some $k,l\geq2$.
\end{thm}

\begin{proof} Let $\{f_k\}_{k\ge1}\subset \Lambda(\psi)$ be a normalized weakly null sequence of
independent identically distributed random variables on $[0,1]$.
Note that we automatically have $\int_0^1f_kdm=0$, $k\ge 1$.

Using standard symmetrization trick, we consider another sequence
$\{f_k'\}_{k\ge 1}$ of independent random variables (which is also
independent with respect to the sequence $\{f_k\}_{k\ge 1}$) such
that $f_k'$ is equidistributed with $f_k$ and define
$h_k:=f_k-f_k'$, $k\ge 1$. Clearly, $\{h_k\}_{k\ge 1}$ is a
sequence of independent symmetric identically distributed random
variables. Noting, that by \cite[Proposition~11, p.~6]{Br1994}, we
have
$$||\sum_{k=1}^nf_k||_{\Lambda(\psi)}\leq const\cdot
||\sum_{k=1}^nh_k||_{\Lambda(\psi)},\quad n\ge 1.$$ Now, if $\psi$
satisfies conditions \eqref{first property} and \eqref{second
property}, then it follows from Corollary \ref{lorentz
alternative} that $||\sum_{k=1}^nh_k||_{\Lambda(\psi)}\leq
const\cdot n^q$ for some $q\in (0,1)$ and hence $\frac{1}{q}\in
\Gamma_{\rm iid}(\Lambda(\psi))$. Conversely, let $\frac{1}{q}\in
\Gamma_{\rm iid}(\Lambda(\psi))$ for some $q\in (0,1)$. Fix $f\in
\Lambda(\psi)$ and consider the sequence $\{(f\otimes r)_n\}_{n\ge
1}\subset \Lambda(\psi)(\Omega, \mathcal{P})$, where the
probability space $(\Omega, \mathcal{P})$ is the infinite direct
product of measure spaces $([0,1],m)$. Since Lorentz spaces
$\Lambda(\psi)(\Omega, \mathcal{P})$ and $\Lambda(\psi)(0,1)$ are
isometric, and since the sequence $\{(f\otimes r)_n\}_{n\ge 1}$ is
weakly null in $\Lambda(\psi)(\Omega, \mathcal{P})$ ( see e.g.
\cite[Lemma 3.4]{SeSu}), we have
$$
\sup_{n\ge 1}\frac{1}{n^q}\|(f\otimes r)_1+(f\otimes r)_2+\dots
+(f\otimes r)_n\|_{\Lambda(\psi)}\leq C(f).
$$
Setting, $B_n:=\frac{1}{n^q}A_n$, $n\ge 1$ we have
$\|B_nf\|_{\Lambda(\psi)}\leq C(f)$ for every $n\ge 1$. By the
uniform boundedness principle, we have
$\|B_n\|_{\mathcal{L}(\Lambda(\psi))}\leq C<\infty$ for all $n\ge
1$, or equivalently that $||A_n||_{\mathcal{L}(\Lambda(\psi))}\leq
C\cdot n^q,\ n\geq 1$. Corollary 9 now yields that the function
$\psi$ satisfies conditions \eqref{first property} and
\eqref{second property}.
\end{proof}

The following Corollary follows immediately from the above
combined with Corollary \ref{coincidence_in_Lor}.

\begin{cor}\label{mainsecond_add} $\Gamma_{\rm i}(\Lambda(\psi))\neq\{1\}$, if and only if the function $\psi\in \Psi$  satisfies conditions \eqref{first
property} and \eqref{second property} for some $k,l\geq2.$
\end{cor}

%
%

We complete this section with the description of
$\Gamma_i(exp(L_p)_0)$, $1\leq p<\infty$.

\begin{thm}\label{Marc} For every $1\leq p\leq2,$ we have
$\Gamma_{\rm iid}(exp(L_p)_0)=\Gamma_{\rm i}(exp(L_p)_0)=[1,2]$.
For every $2\leq p<\infty,$ we have $\Gamma_{\rm
iid}(exp(L_p)_0)=\Gamma_{\rm i}(exp(L_p)_0)= [1,\frac{p}{p-1}].$
\end{thm}

\begin{proof} The first assertion follows from Remark
\ref{connection},  Theorem \ref{explp bound} and Corollary
\ref{coincidence_in_Lor}. The same argument shows that
$\Gamma_i(exp(L_p)_0)\supseteq [1,\frac{p}{p-1}]$ for every $2\leq
p<\infty$. The equality $\Gamma_{\rm i}(exp(L_p)_0)=
[1,\frac{p}{p-1}]$ follows from the fact that the estimate $$\|A_n
\chi _{[0,1]}\|_{exp(L_p)_0}\leq const\cdot n^{1-1/p},\quad n\ge
1$$ is the best possible (see \cite[Theorem 8]{RS} or
\cite[Theorem 15]{D}).
\end{proof}

\section{Concluding Remarks and Examples}

The preceding theorem shows that the set $\Gamma_i(exp(L_p)_0)$ is
non-trivial for all $1\leq p<\infty$, whereas $exp(L_p)$ has the
Kruglov property if and only if $0<p\leq 1$. This result extends
and complements \cite{AsSeSu2007}, where examples of r.i. spaces
$E$ with Kruglov property such that $\Gamma(E)=\{1\}$ and
$\Gamma_{\rm i}(E)\neq \{1\}$ are built. We now present an example
of Lorentz space $\Lambda(\psi)$ such that $\Gamma_{\rm
i}(\Lambda(\psi))\neq\{1\}$ and which does not possess the Kruglov
property.

\begin{ex} Let $\psi\in \Psi$ be given by the condition
$\psi(t):=\frac{1}{\log^{\frac12}(\frac1t)}$, $t\in
[0,e^{-\frac32}]$ and be linear on $ [e^{-\frac32},1]$. The space
$\Lambda(\psi)$ does not have the Kruglov property, however
$\Gamma_{\rm i}(\Lambda(\psi))\neq\{1\}$
\end{ex}

\begin{proof}


Since for every $k,l>1$ we have
$$\lim_{u\rightarrow0}\frac{\psi(ku)}{\psi(u)}=
\lim_{u\rightarrow0}(\frac{\log(u)}{\log(ku)})^{\frac12}=
1<k, \quad\lim_{u\rightarrow0}\frac{\psi(u^l)}{\psi(u)}=
\lim_{u\rightarrow0}(\frac{\log(u)}{\log(u^l)})^{\frac12}=\frac1{l^{\frac12}}<1$$
we see that $\Gamma_{\rm i}(\Lambda(\psi))\neq\{1\}$ by Corollary
9.

By \cite[Theorem 5.1]{AsSeSu2005} a Lorentz space
$\Lambda(\phi)$, $\phi\in \Psi$ has the Kruglov  property if and
only if
$$\sup_{t>0}\frac1{\phi(t)}\sum_{n=1}^{\infty}\phi(\frac{t^n}{n!})<\infty.$$
In our case, for every fixed $t\leq e^{-\frac32}$
$$\sum_{n=1}^{\infty}\psi(\frac{t^n}{n!})=
\sum_{n=1}^{\infty}\frac1{(\log(n!)+n\log(\frac1t))^{\frac12}}=
\sum_{n=1}^{\infty}\frac1{(n\log(n)(1+o(1)))^{1/2}}=\infty.
$$
\end{proof}

\bigskip

\bigskip
\bigskip
\bigskip

\bigskip

\bigskip
\leftline{F. Sulochev}\leftline{School of Mathematics and Statistics}\leftline{University of New South Wales, Kensington NSW 2052}\leftline{Email Address:{\it f.sukochev@unsw.edu.au}}
\leftline {D. Zanin}\leftline {School of Computer Science, Engineering and Mathematics} \leftline {Flinders University,
Bedford Park, SA 5042, Australia}\leftline {Email Address: {\it
zani0005@csem.flinders.edu.au} }

\end{document}